\documentclass[12pt,reqno]{amsart}
\usepackage{amsmath,amssymb,amscd,latexsym,amsthm,mathrsfs,comment}
\usepackage[unicode]{hyperref}
\usepackage{hypbmsec}
\newtheorem{Thm}{Theorem}
\newtheorem{Con}{Conjecture}
\newtheorem{Lem}{Lemma}
\newtheorem{Ans}{Answer}
\newtheorem{Qu}{Question}

\newtheorem{cor}{Corollary}
\theoremstyle{remark}

\begin{document}

\title{Near-primitive roots}
\author{Pieter Moree} 
\address{Max-Planck-Institut f\"ur Mathematik, Vivatsgasse 7, D-53111 Bonn, Germany}
\email{moree@mpim-bonn.mpg.de}
\begin{abstract}
Given an integer $t\ge 1$, a rational number $g$ and a prime $p\equiv 1({\rm mod~}t)$ we say that $g$ is a near-primitive root of
index $t$ if $\nu_p(g)=0$, and $g$ is of order $(p-1)/t$ modulo $p$. In the case $g$ is not minus
a square we compute the density, under the Generalized Riemann Hypothesis (GRH), of such primes explicitly
in the form $\rho(g)A$, with $\rho(g)$ a rational number and $A$ the Artin constant. We follow in this the
approach of Wagstaff, who had dealt earlier with the case where $g$ is not minus a square. The outcome is
in complete agreement with the recent determination of the density using a 
very different, much more algebraic, approach due to Hendrik Lenstra, the author and Peter Stevenhagen. 
\end{abstract}

\maketitle
\section{Introduction}
\noindent Let $g\in \mathbb Q\backslash \{-1,0,1\}$. Let $p$ be a prime. Let $\nu_p(g)$ denote the exponent of
$p$ in the canonical factorization of $g$. If $\nu_p(g)=0$, then we define
$$r_g(p)=[(\mathbb Z/p\mathbb Z)^*:\langle g~{\rm mod~}p\rangle],$$ that is $r_g(p)$ is the residual index
modulo $p$ of $g$. Note that $r_g(p)=1$ iff $g$ is a primitive root modulo $p$. 
For any natural number $t$, let $N_{g,t}$
denote the set of primes $p$ with $\nu_p(g)=0$ and $r_g(p)=t$ (that is $N_{g,t}$ is
the set of near-primitive roots of index $t$). 
Let $\delta(g,t)$ be the natural density of this set of primes (if it exists).
For arbitrary real $x>0$, we
let $N_{g,t}(x)$ denote the number of primes $p$ in $N_{g,t}$ with $p\le x$.\\ 
\indent In 1927 Emil Artin conjectured that for $g$ not equal  to $-1$ or a square, the set
$N_{g,1}$ is infinite and that $N_{g,1}(x)\sim c_gA\pi(x)$, with $c_g$ an explicit rational
number,
$$A=\prod_p\Big(1-{1\over p(p-1)}\Big)\approx 0.3739558,$$
and $\pi(x)$ the number of primes $p\le x$. The constant $A$ is now called Artin's constant.
On the basis of computer experiments by the Lehmers in 1957 Artin had to 
admit that `The machine caught up with me' and provided
a modified version of $c_g$. See e.g.~Stevenhagen \cite{PeterS} for some of
the historical details. On GRH this modified version was shown to be correct by Hooley
\cite{Hooley}.\\ 
\indent Thus $\delta(g,1)$ is explicitly known (under GRH). Determining similarly $\delta(g,t)$
turns out to be rather more difficult and for ease of exposition we first consider the case
where $g>1$ is square free. In this case work of Lenstra \cite{Lenstra} and Murata \cite{Murata}
suggests the following conjecture (with as usual $\mu$ the
M\"obius function and $\zeta_k=e^{2\pi i/k}$).
\begin{Con}
\label{leo}
Let $g>1$ be a square free integer and $t\ge 1$ an integer. The set 
$N_{g,t}$ has a natural density $\delta(g,t)$ 
which is given in Table 1. 
We have
$$N_{g,t}{\it ~is~finite~iff~}\delta(g,t)=0 {\rm ~iff~}g\equiv 1({\rm mod~}4),~2\nmid t,~g|t.$$
\end{Con}
We note that if $g\equiv 1({\rm mod~}4)$, $2\nmid t$ and $g|t$, then $N_{g,t}$ is finite.
To see this note that in this case we have $({g\over p})=1$ for the primes
$p\equiv 1({\rm mod~}t)$ by the law of quadratic reciprocity and thus $r_g(p)$ must
be even, contradicting the assumption $2\nmid t$.\\
\noindent Note that if a set of primes is finite, then its natural density is zero. The converse
is often false, but for a wide class of Artin type problems (including the one
under consideration in this note) is true (on GRH) as first
pointed out by Lenstra \cite{Lenstra}.\\
\indent Given an integer $a$ and a prime $q$, we write $a_q$ to denote the $q$-part of $a$ (that is
$a_q=q^{\beta}$ with $q^{\beta}|a$ and $q^{\beta+1}\nmid a$). 
We put 
\begin{equation}
\label{eetee}
B(g,t)=\prod_{p|g,~p\nmid t}{-1\over p^2-p-1},~~
E(t)={A\over t^2}\prod_{p|t}{p^2-1\over p^2-p-1}.
\end{equation}
Note that if $g|t$, then in the definition of $B(g,t)$ we
have the empty product and hence $B(g,t)=1$. It follows that if further $t$ is odd and $g\equiv 1({\rm mod~}4)$, 
then $\delta(g,t)=0$. The maximal value of $\delta(g,t)$ that occurs is $2E(t)$. Table 1 we took from
a paper by Murata \cite{Murata}. We will show that the densities in Table 1 can be compressed into one equation, namely (\ref{oneline}).
\vskip .8mm

\begin{table}[ht]
\centering
{\bf Table 1: The density $\delta(g,t)$ of $N_{g,t}$ (on GRH)}
\vskip 3mm

\begin{tabular}{|r|c|r|}
\hline
$g$ & $t_2$ & 
 $\delta(g,t)$\\ 
\hline
$g\equiv 1({\rm mod~}4)$ & $t_2=1$ & 
$(1-B(g,t))E(t)$\\
  & $2|t_2$ & $(1+B(g,t))E(t)$\\\hline 
$g\equiv 2({\rm mod~}4)$ & $t_2<4$ &  
$E(t)$\\ 
 & $t_2=4$ &  
$(1-B(g,t)/3)E(t)$\\
 & $t_2>4$ &   
$(1+B(g,t))E(t)$\\\hline 
$g\equiv 3({\rm mod~}4)$ & $t_2=1$ &   
$E(t)$\\ 
 & $t_2=2$ &  
$(1-B(g,t)/3)E(t)$\\
 & $t_2\ge 4$ &   
$(1+B(g,t))E(t)$\\\hline 
\end{tabular}
\end{table}
\begin{Thm}
\label{uno}
Conjecture {\rm \ref{leo}} holds true on GRH.
\end{Thm}
The proof is postponed untill Section \ref{sectietwee}.

\section{Generalization to rational $g$}
\label{sectietwee}
\noindent A natural next problem is to study what
happens if one relaxes the condition that $g$ should be square free. 
Our starting point here will be a result due to Wagstaff \cite{W}. We need
some notation. We put
$$S(h,t,m)=\sum_{n=1\atop m|nt}^{\infty}{\mu(n)(nt,h)\over nt\varphi(nt)},$$
with $\varphi$ Euler's totient function.
\begin{Thm} {\rm \cite{W}}. {\rm (GRH)}.
\label{waggel}
\label{leno}
Let $g\in \mathbb Q\backslash \{-1,0,1\}$ and $t\ge 1$ be an arbitrary integer. Write
$g=\pm g_0^h$, where $g_0\in \mathbb Q$ is positive and not an exact power of a rational
and $h\ge 1$ an integer. Let $d(g_0)$ denote the discriminant of $\mathbb Q(\sqrt{g_0})$.
The natural density of the set $N_{g,t}$, $\delta(g,t)$, exists 
and is given by
\begin{equation}
\label{musum}
\sum_{n=1}^{\infty}{\mu(n)\over [\mathbb Q(\zeta_{nt},g^{1/nt}):\mathbb Q]},
\end{equation}
which
equals a rational number times the Artin constant $A$. 
Write
$g_0=g_1g_2^2$, where $g_1$ is a square free integer and $g_2$ is a rational. If
$g>0$, set $m={\rm lcm}(2h_2, d(g_0))$. For $g<0$, define $m=2g_1$ if $2\nmid h$
and $g_1\equiv 3({\rm mod~}4)$, or $h_2=2$ and $g_1\equiv 2({\rm mod~}4)$; let
$m={\rm lcm}(4h_2,d(g_0))$ otherwise. If $g>0$, we have
$\delta(g,t)=S(h,t,1)+S(h,t,m)$. If $g<0$ we have
\begin{equation}
\label{vierergruppe}
\delta(g,t)=S(h,t,1)-{1\over 2}S(h,t,2)+{1\over 2}S(h,t,2h_2)+S(h,t,m).
\end{equation}
\end{Thm}
In case $g>0$ or $2\nmid h$, Wagstaff expressed $\delta(g,t)$ as an Euler product. By the
work of Lenstra \cite{Lenstra} we know this is also possible in general. The next theorem
achieves this. Partial inspiration for it came from recent joint work with Lenstra and
Stevenhagen, see Section \ref{prospect}.
\begin{Thm} {\rm (GRH)}.
\label{mainz}
Let $g\in \mathbb Q\backslash \{-1,0,1\}$ and $t\ge 1$ be an arbitrary integer. Write
$g=\pm g_0^h$, where $g_0\in \mathbb Q$ is positive and not an exact power of a rational
and $h\ge 1$ an integer. Let $d(g_0)$ denote the discriminant of $\mathbb Q(\sqrt{g_0})$. Put
$F_p=\mathbb Q(\zeta_p,g^{1/p})$.
Put 
$$A(g,t)={(t,h)\over t^2}\prod_{p|t,~h_p|t_p}(1+{1\over p})\prod_{p\nmid t}(1-{1\over [F_p:\mathbb Q]}).$$
Put
$$\Pi_1=\prod_{p|d(g_0),~p\nmid 2t}{-1\over [F_p:\mathbb Q]-1}.$$
Put 
\begin{equation}
\label{e2m2}
E_2(m_2)=\begin{cases}  
1  & {\rm ~if~}m_2| t_2;\cr
-1/3 & {\rm ~if~}m_2 = 2t_2 \ne 2;\cr
-1  & {\rm ~if~}m_2 = 2t_2 = 2;\cr
0  & {\rm ~if~}m_2 \nmid 2t_2,\cr
\end{cases}
\end{equation}
We have
\begin{equation}
\label{adoora}
{A(g,t)\over A}={(t,h)\over t^2}\prod_{p|th}{1\over p^2-p-1}\prod_{p|t\atop pt_p|h_p}p(p-1)\prod_{p|t\atop h_p|t_p}(p^2-1)\prod_{p|h,~p\nmid t_1}p(p-2),
\end{equation}
where
$$t_1=\begin{cases}
2t & {\rm ~if~}g<0,~2|h,~2\nmid t;\cr
t & {\rm ~otherwise}.
\end{cases}$$
Note that $A(g,t)=0$ iff $g>0$, $2|h$ and $2\nmid t$.\\
\indent The natural density of the set $N_{g,t}$ exists, denote it by $\delta(g,t)$.\\
Put $v_0={\rm lcm}(2h_2,d(g_0)_2)$ and $v={\rm lcm}(2h_2,d(g)_2)$.\\
If $g>0$, then $\delta(g,t)=A(g,t)(1+E_2(v_0)\Pi_1)$.\\ If $h$ is odd, then
$\delta(g,t)=A(g,t)(1+E_2(v)\Pi_1)$.\\
If $g<0$, $2|h$ and $2\nmid t$, we have
$\delta(g,t)=A(g,t)$.\\
\indent Next assume $g<0$, $2|(h,t)$.\\
If $h_2=2$ and $8|d(g_0)$, then
\begin{equation}
\label{dubbel} 
\delta(g,t)=\begin{cases}
{1\over 3}A(g,t)(1-\Pi_1) & {\rm ~if~}t_2=2;\cr
A(g,t)(1+\Pi_1) & {\rm ~if~}4|t_2.\cr
\end{cases}
\end{equation}
In the remaining cases we have
$$\delta(g,t)=\begin{cases}
A(g,t)/2 & {\rm ~if~}2t_2|h_2;\cr
A(g,t)/3 & {\rm ~if~}t_2=h_2;\cr
A(g,t)(1-{1\over 3}\Pi_1) & {\rm ~if~}t_2=2h_2;\cr
A(g,t)(1+\Pi_1) & {\rm ~if~}4h_2|t_2.\cr
\end{cases}$$
\end{Thm}
\begin{cor} {\rm (GRH)}. \label{blubbo}
Let $g>1$ be a square free integer. Then
\begin{equation}
\label{oneline}
\delta(g,t)=(1+E_2({\rm lcm}(2,d(g)_2))B(g,t))E(t).
\end{equation}
\end{cor}
{\it Proof}. We have $A(g,t)=S(1,t,1)=E(t)$ (see the remark
following Lemma \ref{t1t2}). Furthermore, if $2|g$ and $2\nmid t$, then
$\Pi_1=-B(g,t)$ and $\Pi_1=B(g,t)$ otherwise. Since 
$E_2({\rm lcm}(2,d(g)_2))=0$ if $g|2$ and $2\nmid t$, we infer that
$E_2({\rm lcm}(2,d(g)_2))\Pi_1=E_2({\rm lcm}(2,d(g)_2))B(g,t)$. Now invoke
the theorem. \qed

\begin{cor} {\rm (GRH)}.
If $t$ is odd, then
$$\delta(g,t)=A(g,t)(1-{1\over 2}(1-(-1)^{h|d(g)|})\Pi_1).$$
\end{cor}
{\tt Remark}. On putting $t=1$ one obtains the classical result of Hooley \cite{Hooley}.\\

{\it Proof of Theorem} \ref{uno}. On distinguishing cases according to the value of $d(g)_2$, 
Corollary \ref{blubbo} yields Table 1. {}From Table 1 one easily reads off that if $\delta(g,t)=0$,
then $2\nmid t$, $g\equiv 1({\rm mod~}4)$ and $g|t$. In this case we have $(g/p)=1$ for the
primes $p\nmid g$ with $p\equiv 1({\rm mod~}t)$ by the law of quadratic reciprocity and hence
$N_{g,t}$ is finite and so $\delta(g,t)=0$. \qed\\

\noindent The proof of Theorem \ref{mainz} will be given in Section \ref{vier}. It will make use
of properties of Wagstaff sums that will be established in the next section.

\section{Bringing the Wagstaff sums in Euler product form}
Recall the definition of the Wagstaff sum
$$S(h,t,m)=\sum_{n=1\atop m|nt}^{\infty}{\mu(n)(nt,h)\over nt\varphi(nt)}.$$
A trivial observation is that if the divisibility condition forces $n$ to be non-square free, then
$\mu(n)=0$ and hence $S(h,t,m)=0$. This happens for example if $m_2\nmid 2t_2$ (cf. Lemma 
\ref{shtminprod}).

In case $m=1$ it is easily written as an Euler product (here we use that $\mu$ and $\varphi$ are
multiplicative functions).
\begin{Lem}
\label{t1t2}
{\rm 1)}. We have
$$S(h,t,1)={(t,h)\over t^2}\prod_{p|t,~h_p|t_p}(1+{1\over p})\prod_{p\nmid t}\Big(1-{(p,h)\over p(p-1)}\Big).$$
In particular, $S(h,t,1)=0$ iff $2|h$ and $2\nmid t$.\\
{\rm 2)}. If $2|h$ and $2\nmid t$, then
$$S(h,t,2)=-{(t,h)\over t^2}\prod_{p|t,~h_p|t_p}(1+{1\over p})\prod_{p\nmid 2t}\Big(1-{(p,h)\over p(p-1)}\Big).$$
\end{Lem}
{\it Proof}. 1) We have
$$S(h,t,1)={(t,h)\over t\varphi(t)}\sum_n {\mu(n)(nt,h)\varphi(t)\over n\varphi(nt)(t,h)}
={(t,h)\over t\varphi(t)}\prod_p\Big(1-{(pt,h)\varphi(t)\over p\varphi(pt)(t,h)}\Big),$$
where we used that the sum $S(h,t,1)$ is absolutely convergent and the fact that the argument
in the second sum is a multiplicative function in $n$. The contribution of the primes dividing
$t$ to this product is
$${(t,h)\over t\varphi(t)}\prod_{p|t,~pt_p|h_p}(1-{1\over p})\prod_{p|t,~h_p|t_p}(1-{1\over p^2})
={(t,h)\over t^2}\prod_{p|t,~h_p|t_p}(1+{1\over p}),$$
where we used that $\varphi(t)/t=\prod_{p|t}(1-1/p)$.
If $p\nmid t$, then 
$$1-{(pt,h)\varphi(t)\over p\varphi(pt)(t,h)}=1-{(p,h)\over p(p-1)},$$
and part 1 follows.\\
2) We have
$$S(h,t,2)=\sum_{2|n}{\mu(n)(nt,h)\over nt\varphi(nt)}=-\sum_{2\nmid n}{\mu(n)(nt,h)\over nt\varphi(nt)}.$$
The latter sum has the same Euler product as $S(h,t,1)$, but with the factor for $p=2$ omitted.
\qed\\

{\tt Remark}. The above lemma and the definition of the Artin constant shows that $E(t)=S(1,t,1)$ and
$A=S(1,1,1)$.\\

Write $M=m/(m,t)$ and $H=h/(Mt,h)$. Then we have \cite[Lemma 2.1]{W}
$$S(h,t,m)=\mu(M)(Mt,h)E(t)\prod_{q|(M,t)}{1\over q^2-1}\prod_{q|M\atop q\nmid t}{1\over q^2-q-1}
\prod_{q|(t,H)\atop q\nmid M}{q\over q+1}\prod_{q|H\atop q\nmid Mt}{q(q-2)\over q^2-q-1}.$$
The parameter $H$ can be avoided as the formula can be rewritten as
\begin{equation}
\label{langlang}
{\mu(M)(Mt,h)A\over t^2}\prod_{q|mth}{1\over q^2-q-1}\prod_{q|t,~qt_q|h_q\atop m_q|t_q}q(q-1)
\prod_{q|t,~h_q|t_q\atop m_q|t_q}(q^2-1)\prod_{q|h\atop q\nmid mt}q(q-2).
\end{equation}
(In order to see this it is helpful to consider the cases $m_q|t_q$, that is $M_q=1$, and $qt_q|m_q$, that is $q|M$, separately.)
These formulae relate $S(h,t,m)$ to $S(1,t,1)$ ($=E(t)$), respectively to $S(1,1,1)$ ($=A$), however, as we will show, 
expressions simplify considerably if we relate $S(h,t,m)$ to $S(h,t,1)$. We start by showing how to
remove odd prime factors from $m$.
\begin{Lem}
\label{graadoverkop}
Suppose that $p\nmid 2m$.
Then $$S(h,t,mp)=\begin{cases}-S(h,t,m)/({p(p-1)\over (p,h)}-1) & {\rm ~if~}p\nmid t;\cr
S(h,t,m) & {\rm ~if~}p|t.\cr
\end{cases}$$
\end{Lem}
{\it Proof}. If $p|t$ the summation condition $mp|nt$ in the definition of $S(h,t,mp)$ is equivalent
with $m|nt$, that is we have $S(h,t,mp)=S(h,t,m)$.\\ 
\indent Next assume that $p\nmid t$. We have
$$S(h,t,mp)=\sum_{m|nt\atop p|n}{\mu(n)(nt,h)\over nt\varphi(nt)}=\sum_{m|nt}{\mu(pn)(pnt,h)\over pnt\varphi(pnt)}=-{(p,h)\over p(p-1)}\sum_{m|nt\atop p\nmid n}{\mu(n)(nt,h)\over nt\varphi(nt)}.$$
On noting that the latter sum can be written as $S(h,t,m)-S(h,t,mp)$, the proof is then completed.\qed

\begin{Lem}
\label{twooh}
Suppose that we are not in the case where $h$ is even and $t$ is odd.
We have
$$S(h,t,2t_2)=\begin{cases}
-S(h,t,1)/3 & {\rm ~if~lcm}(2,h_2)|t_2;\cr
-S(h,t,1) & {\rm ~if~lcm}(2,h_2)\nmid t_2.
\end{cases}$$
\end{Lem}
{\it Proof}. We can write
$$S(h,t,2t_2)=\sum_{2|n}{\mu(n)(nt,h)\over nt\varphi(nt)}=-{1\over 2}\sum_{2\nmid n}{\mu(n)(2nt,h)\over nt\varphi(2nt)}=\epsilon \sum_{2\nmid n}{\mu(n)(nt,h)\over nt\varphi(nt)},$$
where $\epsilon$ is easily determined (and $\epsilon\ne -1$). Since the latter sum is equal to
$S(h,t,1)-S(h,t,2t_2)$, we then infer that $S(h,t,2t_2)={\epsilon\over 1+\epsilon}S(h,t,1)$. Working
out the remaining details is left to the reader. \qed

\begin{Lem}
\label{shtminprod}
Let $m$ be an integer, having square free odd part. Let $h$ and $t$ be integers, with the requirement
that $t$ be even in case $h$ is even. Then
$$S(h,t,m)=S(h,t,1)E_1(m_2)\prod_{p|m,p\nmid 2t}{-1\over {p(p-1)\over (p,h)}-1},$$
where $$E_1(m_2)=\begin{cases}  
1  & {\rm ~if~}m_2| t_2\cr
-1/3 & {\rm ~if~}m_2 = 2t_2{\rm ~and~lcm}(2,h_2)|t_2\cr
-1  & {\rm ~if~}m_2 = 2t_2{\rm ~and~lcm}(2,h_2)\nmid t_2\cr
0  & {\rm ~if~}m_2\nmid 2t_2,\cr
\end{cases}$$
In case $2h_2|m_2$, we have $E_1(m_2)=E_2(m_2)$, where $E_2(m_2)$ is given
by {\rm (\ref{e2m2})}.
\end{Lem}
{\it Proof}. By Lemma \ref{t1t2} the conditions imposed on $h$ and $t$
imply that $S(h,t,1)\ne 0$. 
By Lemma \ref{graadoverkop} it suffices to show that
$S(h,t,m_2)=S(h,t,1)E_1(m_2)$. If $m_2|t_2$, then no divisibility condition on $n$ is
imposed in the definition of $S(h,t,m_2)$ and so we obtain $S(h,t,m_2)=S(h,t,1)$ and
hence $E_1(m_2)=1$. In case $m_2=2t_2$ we invoke Lemma \ref{twooh}.
If $m_2\nmid 2t_2$, then the summation condition $m|nt$ implies $4|n$ and hence
$\mu(n)=0$ and so $S(h,t,m_2)=0$ and hence $E_1(m_2)=0$.\\
\indent The final claim follows on noting that if $2h_2|m_2$ and $m_2=2t_2$, then
$h_2|t_2$ and hence lcm$(2,h_2)\nmid t_2$ iff $2\nmid t_2$. \qed

\section{Proof of Theorem \ref{mainz}}
\label{vier}
The idea of the proof is to express $\delta(g,t)$ in terms of $S(h,t,1)$, except in case $g<0$, $2|h$ and
$2\nmid t$, when $S(h,t,1)=0$, in which case we express $\delta(g,t)$ in terms of $S(h,t,2)$. 
These two Wagstaff sums are then related to $A(g,t)$ using the following lemma. Note that it shows that
the dependence of $A(g,t)$ on $g$ is weak, as only $h$ and the sign of $g$ matter.
\begin{Lem}
\label{stoa}
We have
$$A(g,t)=\begin{cases}
-S(h,t,2)/2 & {\rm ~if~}g<0,~2|h,~2\nmid t;\cr
S(h,t,1) & {\rm ~otherwise}.\cr
\end{cases}$$
\end{Lem}
{\it Proof}. Note that if $g<0$ and $2|h$, then $F_2=\mathbb Q(i)$ and $[F_2:\mathbb Q]=2$. 
In the remaining cases we have $[F_p:\mathbb Q]=p(p-1)/(p,h)$. On invoking Lemma \ref{t1t2} the
proof is then completed. \qed

\noindent {\it Proof of Theorem} \ref{mainz}. 
Equation (\ref{adoora}) follows by Lemma \ref{stoa} and (\ref{langlang}). We will use a few times, cf. the
proof of Lemma \ref{stoa}, that
$$\Pi_1=\prod_{p|d(g_0),~p\nmid 2t}{-1\over [F_p:\mathbb Q]-1}=\prod_{p|d(g_0),~p\nmid 2t}{-1\over {p(p-1)\over 
(p,h)}-1}.$$
Assume GRH.\\
\noindent {\tt The case $g>0$}.\\ 
By Theorem \ref{waggel} we have
$\delta(g,t)=S(h,t,1)+S(h,t,m)$, with $m={\rm lcm}(2h_2,d(g_0))$.
First assume that $2|h$ and $2\nmid t$. 
Then, by Lemmas \ref{t1t2}  and \ref{stoa}, we have $S(h,t,1)=A(g,t)=0$ and we need
to show that $\delta(g,t)=0$. Since $S(h,t,1)=0$ it remains to show that $S(h,t,m)=0$.
Since for the $n$ in the summation we have $4|2h_2|n$, this is clear. 
Next assume we are in the remaining case, that is either $h$ is odd, or
$2|(h,t)$.  Then $S(h,t,1)=A(g,t)$ by Lemma \ref{stoa}. Note that $m_2=v_0$.
By Lemma \ref{shtminprod} we then find that $\delta(g,t)=S(h,t,1)(1+E_1(v_0)\Pi_1)=A(g,t)(1+E_2(v_0)\Pi_1)$, 
where we have used that $2h_2|v_0$.\\
{\tt The case $h$ is odd}.\\ 
If $g>0$ then $v=v_0$ and we are done, so assume that $g<0$.
The formula for $m$ in Theorem \ref{waggel} can be rewritten as lcm$(2,|d(g)|)$, and one finds that
$\delta(g,t)=S(h,t,1)+S(h,t,{\rm lcm}(2,|d(g)|))$. This is the same
formula as in case $g>0$ and $2\nmid h$, but with $d(g_0)$ replaced by $|d(g)|$. On noting that
the odd part of $d(g_0)$ equals the odd part of $d(g)$, the result then follows.\\
{\tt The case $g<0$, $2\nmid t$ and $2|h$}.\\
We have $S(h,t,1)=S(h,t,m)=S(h,t,2h_2)=0$ and hence $\delta(g,t)=-S(h,t,2)/2$ by (\ref{vierergruppe}). 
Now invoke Lemma \ref{stoa} to obtain $\delta(g,t)=A(g,t)$.\\
{\tt The case $g<0$ and $2|(h,t)$}.\\
Note that $2|m$ and $S(h,t,1)=A(g,t)$. By Lemma \ref{shtminprod} we infer that
$S(h,t,2)=S(h,t,1)$ and $S(h,t,2h_2)=S(h,t,1)E_2(2h_2)$, where
$$E_2(2h_2)=\begin{cases}
1 & {\rm ~if~}2h_2|t_2;\cr
-1/3 & {\rm ~if~}h_2=t_2;\cr
0 & {\rm ~if~}h_2\nmid t_2.\cr
\end{cases}
$$
Note that $$E_2(4)=\begin{cases}
1 & {\rm ~if~}4|t_2;\cr
-1/3 & {\rm ~if~}t_2=2.\cr
\end{cases}
$$
If $h_2=2$ and $8|d(g_0)$, then by Theorem \ref{waggel} we 
have $m=2g_1$, which can be rewritten as $m=d(g_0)/2$ (thus $m_2=4$)
and so
\begin{align}
\delta(g,t) & = S(h,t,1)-{S(h,t,2)\over 2}+{S(h,t,4)\over 2}+S(h,t,{d(g_0)\over 2})\nonumber\cr
 & = S(h,t,1)\Big({1\over 2}+{E_2(4)\over 2}+E_2(4)\Pi_1\Big).\nonumber
\end{align}
where we used that, by Lemma \ref{shtminprod}, $S(h,t,d(g_0)/2)=S(h,t,1)E_2(4)\Pi_1$. Using 
that $S(h,t,1)=A(g,t)$ and the formula
for $E_2(4)$, we then arrive at (\ref{dubbel}).\\
\indent In the remaining case, 
$m={\rm lcm}(4h_2,d(g_0))$.
Note that $m_2=4h_2$ and
$$E_2(4h_2)=\begin{cases}
1 & {\rm ~if~}4h_2|t_2;\cr
-1/3 & {\rm ~if~}2h_2=t_2;\cr
0 & {\rm ~if~}2h_2\nmid t_2.
\end{cases}
$$
We find that
$$
\begin{aligned}
\delta(g,t) & =S(h,t,1)-{S(h,t,2)\over 2}+{S(h,t,2h_2)\over 2}+S(h,t,m)\nonumber\cr
& =S(h,t,1)\Big({1\over 2}+{E_2(2h_2)\over 2}+
E_2(4h_2) \Pi_1\Big).\nonumber\cr
\end{aligned}$$
Using that $S(h,t,1)=A(g,t)$ and the formulae for $E_2(2h_2)$ and $E_2(4h_2)$ given above, the proof
is then completed. \qed

\section{Vanishing of $\delta(g,t)$}
The aim of this section is to give a new proof of Theorem \ref{sixcondi} (due to Lenstra \cite{Lenstra}, who stated it without proof). The
first published proof was given by Moree in \cite{near1}.
He introduced a function $w_{g,t}(p)\in \{0,1,2\}$ for which he proved 
(see \cite{near1}, for a rather easier reproof see \cite{near2}) under GRH
that 
$$N_{g,t}(x)=(h,t)\sum_{p\le x,~p\equiv 1({\rm mod~}t)}w_{g,t}(p){\varphi((p-1)/t)\over p-1}
+O\Big({x\log \log x\over \log^2 x}\Big).$$
This function $w_{g,t}(p)$ has the property that, under GRH, $w_{g,t}(p)=0$ for
all primes $p$ sufficiently large iff $N_{g,t}$ is finite. Since the definition of
$w_{g,t}(p)$ involves nothing more than the Legendre symbol, it is then not difficult
to arrive at the cases 1-6. E.g. in case 1 $g$ is a square
modulo $p$, and thus $2|t$, contradicting $2\nmid t$. Likewise for the other 5 cases
the obstructions can be written down (it turns out $r_g(p)_2\ne t_2$ in each case). For the complete list 
of obstructions we refer to
Moree \cite[pp. 170-171]{near1}.\\ 
\indent Regarding the 
six vanishing cases Wagstaff \cite[p. 143]{W} wrote: `It is easy to verify directly that our expression for $\delta(g,t)$
vanishes in each of Lenstra's cases, but it is tedious to check that these are the only cases in which
it vanishes'. We will show that once Wagstaff's result is brought into Euler product form, as
done in Theorem \ref{mainz}, it is straightforward to establish Theorem \ref{sixcondi}. A more
conceptual, shorter and elegant (but less elementary) proof of Theorem \ref{sixcondi} will appear in \cite{LMS}. 
\begin{Thm} {\rm (GRH)}. 
\label{sixcondi}
The set $N_{g,t}$ is finite  iff $\delta(g,t)=0$ iff we are in one of the following six 
(mutually exclusive) cases:\\
{\rm 1)} $2\nmid t$,~$d(g)|t$.\\
{\rm 2)} $g>0$, $2h_2|t_2$, $3\nmid t$, $3|h$, $d(-3g_0)|t$.\\
{\rm 3)} $g<0$, $h_2=1$, $t_2=2$, $3\nmid t$, $3|h$, $d(3g_0)|t$.\\
{\rm 4)} $g<0$, $h_2=2$, $t_2=2$, $d(2g_0)|2t$.\\
{\rm 5)} $g<0$, $h_2=2$, $t_2=4$, $3\nmid t$, $3|h$, $d(-6g_0)|t$.\\
{\rm 6)} $g<0$, $4h_2|t_2$, $3\nmid t$, $3|h$, $d(-3g_0)|t$.
\end{Thm}
\noindent {\tt Example}. (GRH). If $g>1$ is square free, then case 1 is the only one to take
into account and we find $\delta(g,t)=0$ iff $2\nmid t$, $d(g)|t$, that is iff $2\nmid t$, $g|t$, 
$g\equiv 1({\rm mod~}4)$.

\begin{table}[ht]
\centering
{\bf Table 2: Examples of pairs $(g,t)$ satisfying cases 1-6}

\begin{tabular}{|r|c|r|r|c|c|c|}
\hline
& 1 & 2 & 3 & 4 & 5 & 6 \\ \hline
$(g,t)$ & $(5,5)$ & $(3^3,4)$ &  $(-15^3,10)$ & $(-6^2,6)$ & $(-6^6,4)$& $(-3^3,4)$\\ \hline
\end{tabular}
\end{table}

{\it Proof of Theorem} \ref{sixcondi}. If one of 1-6 is satisfied, then
$N_{g,t}$ is finite. This can be shown by elementary
arguments only involving quadratic reciprocity (see Moree \cite[pp. 170-171]{near1}). It is thus enough to show that $\delta(g,t)=0$
iff one of the six cases is satisfied. 
For the proof we will split up case 6 into two subcases:\\
6a) $g<0$, $2|h_2$,~$4h_2|t_2$, $3\nmid t$, $3|h$, $d(3g_0)|t$.\\
6b) $g<0$, $h_2=1$,~$4|t_2$, $3\nmid t$, $3|h$, $d(3g_0)|t$.\\
(For our proof it is more natural to require $d(3g_0)|t$, which, since $4|t$, is equivalent
with $d(-3g_0)|t$.)
Let us denote by $d^{*}(g_0)$ the odd part of the discriminant of $g_0$, that
is $d^{*}(g_0)=d(g_0)/d(g_0)_2$. Note that
\begin{equation}
\label{pie1}
\Pi_1=\begin{cases}
1 & {\rm ~if~}d^{*}(g_0)|t;\cr
-1 & {\rm ~if~}3|d(g_0),~d^{*}(g_0)|3t,~3\nmid t,~3|h;\cr
\in (-1,1) & {\rm ~otherwise}.\cr
\end{cases}
\end{equation}
{\tt The case $2\nmid t$}.\\
If $2|h$ one has $\delta(g,t)=0$ iff $g>0$, that is iff $d(g)|t$.\\
If $2\nmid h$, then $A(g,t)\ne 0$ and we have $\delta(g,t)=0$ iff
$E_2({\rm lcm}(2,d(g)_2))=-1$ and $\Pi_1=1$, that is iff lcm$(2,d(g)_2)=2$
and $d^{*}(g)|t$, that is iff $d(g)|t$.\\
Thus from now on we may assume that $2|t$. This ensures that $A(g,t)\ne 0$.\\
{\tt The case $g>0$ and $2|t$}.\\ 
Now the possibility $E_2(m_2)=-1$ cannot occur and
thus $\delta(g,t)=0$ iff $E_2(m_2)=1$ and $\Pi_1=-1$. The latter two conditions are
both satisfied iff lcm$(2h_2,d(g_0)_2)|t_2$, $3|d(g_0)$, $d^{*}(g_0)|3t$, $3\nmid t$,~$3|h$.
These conditions can be reformulated as $2h_2|t_2$, $3|d(g_0)$,~$d(g_0)|3t$, $3\nmid t$ and $3|h$.
Since $3\nmid t$,~$3|d(g_0)$,~$d(g_0)|3t$ iff $d(-3g_0)|t$, $3\nmid t$, we are done.\\
\indent Thus if $g>0$ or $2\nmid t$, then $\delta(g,t)=0$ iff we are in case 1 or in case 2.
It remains to consider the case where $g<0$ and $2|t$.\\
{\tt The case $g<0$, $2|t$, $2\nmid h$}.\\
Here we have $\delta(g,t)=0$ iff $E_2(v)=1$ and $\Pi_1=-1$.
Note that $E_2(v)=1$ means that we require lcm$(2,d(g)_2)|t_2$.\\
If $t_2=2$, then lcm$(2,d(g)_2)|t_2$ and
$\Pi_1=-1$ iff we are in case 3.\\
If $4|t_2$, then lcm$(2,d(g)_2)|t_2$ and
$\Pi_1=-1$ iff we are in case 6b.\\
{\tt The case $g<0$, $2|(h,t)$}.\\
We have $\delta(g,t)=0$ iff we are in one of the following three cases:\\
A) $h_2=2,~t_2=2,~8|d(g_0),~\Pi_1=1$;\\
B) $h_2=2,~t_2=4,~8|d(g_0),~\Pi_1=-1$;\\
C) $2|h_2,~4h_2|t_2,~\Pi_1=-1$.\\
It is easily checked that these are merely cases 4, 5 and 6a in different guises. \\
\indent To sum up, we have shown that $\delta(g,t)=0$ iff we are in one of the cases
1,2,3,4,5,6a or 6b. Note that the six cases are mutually exclusive.\qed\\

We now propose a conjecture on $\delta(g,t)$ for arbitrary rational $g$. It generalizes
Conjecture 1.
\begin{Con}
\label{netvorige}
The set $N_{g,t}$ has a natural density $\delta(g,t)$ that is given as in Theorem \ref{mainz} and
is a rational multiple of the Artin constant $A$.
The set $N_{g,t}$ is finite iff $\delta(g,t)=0$ iff we are in one of the six cases of 
Theorem \ref{sixcondi}.
\end{Con}
On combining Theorem \ref{mainz} and Theorem \ref{sixcondi} we deduce that Conjecture \ref{netvorige} holds true
on GRH.
\begin{Thm}
Conjecture \ref{netvorige} is true under GRH.
\end{Thm}

\section{Near-primitive roots density through character sum averages}
\label{prospect}
Lenstra, Moree and Stevenhagen \cite{LMS} show that for a large class of Artin-type problems the
set of primes has a natural density $\delta$ that is given by
\begin{equation}
\label{pete}
\delta=(1+\prod_p E_p)\prod_p A_p,
\end{equation}
where $\prod_p A_p$ is the `generic answer' to the density problem (e.g. $A$ in the original Artin
problem) and $1+\prod_p E_p$ a correction factor. For finitely many primes $p$ one has $E_p\ne 1$ and
further $-1\le E_p\le 1$ as $E_p$ is a (real) character sum average over a finite set (and hence
the correction factor is a rational number). In particular, it is rather easy in this set-up to
determine when $\delta=0$. The character sum method makes use of the theory of radical entanglement as 
developped by Lenstra \cite{SanAnton}\\
\indent For the near-primitive root problem the method leads rather immediately to the formula
$\delta(g,t)=A(g,t)(1+E_2'\Pi_1)$ in case $g>0$. The only harder part is the determination of $E_2'$. For
the details the reader is referred to \cite{LMS}.\\
\indent Indeed, the great advance of the newer method is that it very directly leads to a formula for
the density in Euler product form. The classical method leads to infinite sums involving the
M\"obius function and nearly multiplicative functions (in our case Wagstaff's result (Theorem \ref{waggel}). It 
then requires rather cumbersome manipulations to arrive at a density in Euler product form. Indeed, inspired by the predicted result (\ref{pete}) the
author attempted (and managed) to bring Wagstaff's result in Euler product form.\\
\indent The analogue of Theorem \ref{mainz} obtained in this approach, Theorem 6.4 of
\cite{LMS}, looks slightly different from Theorem \ref{mainz}. However, on noting that
$s_2$ as defined in Theorem 6.4 is merely the 2-part of  $m$ as defined in Wagstaff's result
Theorem \ref{waggel}, it is not
difficult to show that both methods give rise to the same Euler products for the density. By allowing
$g_0$ to be negative in case $h$ is odd and $g<0$, the above 6 cases where vanishing occurs can be reduced
to 5 cases (see Corollary 6.5 of \cite{LMS}).

\section{An application}
\noindent Let $\Phi_n(x)$ denote the $n$-th cyclotomic polynomial. Let $S$ be the set of primes $p$ such that
if $f(x)$ is any irreducible factor of $\Phi_p(x)$ over $\mathbb F_2$, then $f(x)$ does not divide
any trinomial. Over $\mathbb F_2$, $\Phi_p(x)$ factors into $r_2(p)$ irreducible polynomials.
Let
$$S_1=(\{p>2:2\nmid r_2(p)\}\}\cup \{p>2:2\le r_2(p)\le 16\})\backslash \{3,7,31,73\}.$$
\begin{Thm}
We have $S_1\subseteq S$.
The set $S_1$ contains the primes $p>3$ such that $p\equiv \pm 3({\rm mod~}8)$.
On GRH the set $S_!$ has density 
\begin{equation}
\label{densiedens}
\delta(S_1)={1\over 2}+A{1323100229\over 1099324800}\approx 0.950077195\cdots 
\end{equation}
\end{Thm}
\noindent {\it Proof}. The set $\{p>2:2\nmid r_2(p)\}\}$ equals the set of primes $p$ such 
that $({2\over p})=-1$, that is the set of primes $p$ such that $p\equiv \pm 3({\rm mod~}8)$. This set has density $1/2$. We thus
find, on consulting Table 1, that
\begin{eqnarray}
\delta(S_1) & = &{1\over 2}+\sum_{2\le j\le 16\atop 2|j}A(2,j)\cr
& = &{1\over 2}+E(2)(1+{2\over 3\cdot 4}+{2\over 16}+{2\over 64})+E(6)(1+{2\over 3\cdot 4})+E(10)+E(14),\nonumber
\end{eqnarray}
which yields (\ref{densiedens}) on invoking the definition (\ref{eetee}) of $E(t)$.
That $S_1\subseteq S$ is a consequence of the work of Golomb and Lee \cite{GL}. \qed\\

\noindent {\tt Acknowledgment}. Given expressions like (\ref{langlang}), my intuition was that expressing
$\delta(g,t)$ in Euler product form would lead to very unpleasant formulae and thus I never attempted this. Discussions with
Peter Stevenhagen, considering the near-primitive root problem by a much more algebraic method, strongly 
suggested easier expressions for $\delta(g,t)$ than expected. 
This led me to try to bring Wagstaff's result in Euler product form, also with the aim of verifying
the results found by the character sum method (alluded to in Section \ref{prospect}).\\
\indent The memorable discussions
with Stevenhagen are gratefully acknowledged. Finally, I like to thank Carl Pomerance for helpful
e-mail correspondence.

\end{document}